\documentclass[11pt]{amsart}
\usepackage{a4wide}
\usepackage{amssymb}
\usepackage{mathrsfs}
\theoremstyle{plain}
\newtheorem{theorem}{Theorem}
\newtheorem*{theorem*}{Theorem}%
\newtheorem{lemma}{Lemma}[section]
\newtheorem{proposition}{Proposition}[section]

\newtheorem{definition} {Definition} [section]
\newtheorem{example} {Example}[section]
\newcommand{\NN}{\mathbf{N}}

\newcommand{\FF}{\mathbf{F}}
\newcommand{\Gd}{\mathbf{G}_d}
\newcommand{\cL}{\mathcal{L}} \newcommand{\cK}{\mathcal{K}}
\newcommand{\cG}{\mathcal{G}} \newcommand{\cH}{\mathcal{H}}
\newcommand{\cP}{\mathcal{P}} 
\newcommand{\urd} {U^r_d}
\newcommand{\e} {\epsilon}\newcommand{\cT} {\mathcal{T} }
\newcommand{\T} {\mathcal{T} }
\newcommand{\SSS} {\mathcal{S}}
\newcommand{\grd}{\mathbf{Gr}_d}
\newcommand{\hg} {\hat{G}}
\newcommand{\urdb}{U^{r,B}_d}
\newcommand{\gdb} {\mathbf{G}^B_d}
\newcommand{\limn}{\lim_{n\to\infty}}
\newcommand{\grdb}{\mathbf{Gr}^B_d}
\newcommand{\grdbk}{\mathbf{Gr}^{B^2}_d}

\newcommand{\urkd} {U^{r,k}_d}

\begin{document}

\title[]{Finite graphs and amenability}
\author[G. Elek]{G\' abor Elek}
\thanks{Work supported in part by a Marie Curie grant and TAMOP 4.2.1/B-09/1/KMR-2010-003}
\begin{abstract} Hyperfiniteness or amenability of measurable equivalence 
relations and group actions has been studied for almost fifty years. Recently,
unexpected applications of hyperfiniteness were found in computer science 
in the
context of testability of graph properties. In this paper we propose a unified 
approach to 
hyperfiniteness. We establish some new results and 
give new proofs of theorems
of Schramm, Lov\'asz, Newman-Sohler and Ornstein-Weiss.

\end{abstract}
\maketitle
\tableofcontents
\section{Introduction}

\subsection{Local statistics for graphs and graphings}
First let us recall some basic notions. 
Let $\Gd$ denote the set of finite simple graphs of vertex degree bound 
$d$ (up
to isomorphism). A rooted graph $H$ of radius at most $r$ is
\begin{itemize}
\item a graph with vertex degree bound $d$ and a distinguished vertex $x$ (the
  root)
\item such that $d_G(x,y)\leq r$ for any $y\in V(G)$, where $d_G$ is the usual
  shortest path metric.
\end{itemize}
Let us denote by $\urd$ the set of all rooted graphs of radius at most $r$ up to
rooted isomorphisms. If $G\in \Gd$ and $\alpha\in \urd$ then
$T(G,\alpha)$ is defined as 
$$T(G,\alpha):=\{v\in V(G)\,\mid\, B_r(v)\sim \alpha\}\,.$$
Set $p(G,\alpha):=\frac{|T(G,\alpha)|}{|V(G)|}\,.$ That is $p(G,\alpha)$
is the probability that the $r$-ball around a random vertex of $G$ is
rooted-isomorphic to $\alpha$.
Let us enumerate the elements of the set
$\cup^\infty_{r=1} \urd$.
Then we get a map $\cL:\Gd\to [0,1]^{\NN}$. We equip $[0,1]^{\NN}$ with
a metric $d_\pi$ that generates the usual product topology.
$$\cL(G)=\{p(G,\alpha_1), p(G,\alpha_2),\dots\} $$
The map is ``almost'' injective: if $\cL(G)=\cL(H)$ then
there exists a graph $K$ such that both $G$ and $H$ are disjoint union 
of $K$-copies.
We say that a sequence of graphs $\{G_n\}\subset \Gd$ is {\it convergent}
(in the sense of Benjamini and Schramm) if $\lim_{n\to\infty} p(G_n,\alpha)$
exists for any $r\geq 1$ and $\alpha\in\urd$. That is $\{G_n\}^\infty_{n=1}$
is convergent if and only if $\{\cL(G_n)\}^\infty_{n=1}$ is convergent 
pointwise.

\noindent
Now we recall the notion of a graphing \cite{Kech}.
Let $X$ be a standard Borel set. A Borel set $E\subset X\times X$ is a {\it
  Borel graph} if
\begin{itemize}
\item $(x,y)\in E$ implies that $(y,x)\in E$
\item $(x,x)\notin E$ if $x\in X$.
\end{itemize}
Note that the degree of a vertex $x$ is well-defined. A Borel graph of vertex
degree bound $d$ is such a graph that all of its components are countable graphs
with vertex degree bound $d$. A measurable graph (or a {\it graphing}) is a
Borel graph on a standard Borel probability measure space $(X,\mu)$ satisfying
the following property.
\begin{itemize}
\item If $T:X\to X$ is a Borel bijection such that either $T(x)=x$ or
$(x,T(x))\in E$, then $T$ preserves the measure $\mu$.
\end{itemize}

 The most important
examples of such graphings are given by group actions. Let $\Gamma$ be a
finitely generated group with a symmetric generating system $S$. Consider a
measure preserving Borel action of $\Gamma$ on $(X,\mu)$. Now let
$(x,y)\in E$ if $x\neq y$ and $sx=y$ for some $s\in S$. Then $\cG(X,E,\mu)$ is
a graphing. For such a graphing $\cG$ with vertex degree bound $d$ we can
define the probabilities $p(\cG,\alpha)$ as well. Let $\alpha\in U^r_d$, then 
$T(\cG,\alpha)$ is the Borel set of points $x\in X$ such that
$B_r(x)\sim\alpha$. Let $p(\cG,\alpha):=\mu(T(\cG,\alpha))\,.$
  Thus we can extend $\cL$ to the isomorphism
 classes of graphings of vertex degree bound $d$ (from now on, all the 
graphing in the paper are  supposed to have vertex degree bound $d$). 
We say that $\cG$
is a limit of a convergent graph sequence $\{G_n\}^\infty_{n=1}\subset \Gd$
if for any $r\geq 1$ and $\alpha\in \urd$
$$\lim_{n\to\infty} p(G_n,\alpha)=p(\cG,\alpha)\,,$$
that is $\lim_n \cL(G_n)=\cL(\cG)$.
We define the pseudo-distance of graphings $d_{stat}(\cG,\cH)$ by
$d_{\pi}(\cL(\cG),\cL(\cH))$.

For any convergent graph sequence
there exists a limit graphing \cite{Elekfin}, the converse statement is
an open conjecture due to Aldous and Lyons \cite{AL}.

\noindent
Let $\cG(X,E,\mu)$ be a graphing and $Z\subset E$ be a Borel set of edges.
Let $$\deg_Z(x):=|\{y\in X\,\mid\,(x,y)\in Z\}|\,.$$
Then $$\mu_E(Z):=\frac{1}{2}\int_X deg_Z(x)\,.$$

\subsection{Hyperfiniteness} \label{hyperke}
The notion of hyperfiniteness was introduced in \cite{Elekspec}.
A set of graphs $\{G_n\}\subset \Gd$ is called a {\it hyperfinite} family if
\begin{itemize}
\item for any $\e>0$ there exists $K>0$ such for each $n\geq 1$ there exists
a set $Z_n\subset V(G_n)$, $|Z_n|<\e|V(G_n)|$ such that if we 
remove the edges incident to $Z_n$
the resulting graph $G'_n$ consists of
components of size at most $K$.
\end{itemize}
Note that any planar or subexponentially growing family of graphs is
hyperfinite \cite{ElekRSA}. Also, F\o lner sequences of a 
finitely generated amenable group form a hyperfinite family. 
Hyperfiniteness can be defined for graphings as well \cite{Kech}.
We call a graphing $\cG$ hyperfinite (or amenable) if for any $\e>0$
there exists $K>0$ such that for some Borel set $Z\subset X$
\begin{itemize}
\item $\mu(Z)<\e$
\item all the components of $E\backslash Z$ have size at most $K$.
\end{itemize}
Note that $E\backslash Z$ denotes the graphing with vertex set $X$, with
edges of $\cG$ that are not incident to an element of $Z$.
The  classical examples of hyperfinite graphings are graphings
of subexponential growth and the ones associated to probability measure
preserving actions of finitely generated amenable groups.
Now we can formulate our first result.
\begin{theorem}\label{elsotetel}
A convergent graph sequence $\{G_n\}^\infty_{n=1}$ is hyperfinite
if and only if its limit graphing $\cG$ is hyperfinite.
\end{theorem}
The original version of this theorem was proved by Oded Schramm \cite{Schramm}
using an ingenious probabilistic idea. Notice that he considered 
{\it unimodular measures} as limit objects. He noted that there is a 
minor technical 
difficulty in some cases (due to symmetries). Our approach is completely 
deterministic and seems to avoid these difficulties.
Interestingly, in both proofs one of the 
directions are much easier to prove than the other, but not the same ones 
(the reason of this strange phenomenon is hidden in the definition of 
hyperfiniteness
for unimodular measures).

\subsection{Equivalences of graphings}
Following Lov\'asz \cite{Lovasz}, we say that two graphings $\cG$ and $\cH$
are weakly equivalent if they have the same local statistics 
$\cL(\cG)=\cL(\cH)$. We will prove the following statement (also proved by
Lov\'asz using Schramm's probabilistic method).
\begin{theorem} \label{masodiktetel}
If $\cG$ and $\cH$ are weakly equivalent then $\cH$ is hyperfinite if and
only if $\cG$ is hyperfinite. That is hyperfiniteness is a local property.
\end{theorem}
We say that two graphings $\cG(X,\mu)$ and $\cH(Y,\nu)$ are strongly 
equivalent if
for any $\epsilon>0$ there exists a measure preserving bijective
map $\rho_\e:X\to Y$ such that 
$$\mu_E(\rho_\e^{-1}E(\cH) \triangle E(\cG))<\e $$ 
If two graphings are strongly equivalent then they are clearly weakly 
equivalent as well. However for hyperfinite graphings, the converse is
true.
\begin{theorem} \label{harmadiktetel}
If $\cG$ and $\cH$ are weakly equivalent hyperfinite graphings, then
they are strongly equivalent.
\end{theorem}
We shall prove a variant of this theorem for group actions as well,
generalizing the classical Rokhlin Lemma.
\subsection{The Equipartition Theorem and its consequences}
The following result states that for a hyperfinite family
statistically similar graphs can be partitioned similarly. 
\begin{theorem}\label{equi}
Let $\cP\subset\Gd$ be a hyperfinite family. Then for any $\e>0$, there exists
an integer $K>0$ with the following property: For any $\delta>0$, 
there exists $f(\delta)>0$ such that
if $G\in\cP$ and $H\in G_d$ with $d_{stat}(G,H)\leq f(\delta)$ then
one can remove less that $2\e |E(G)|$ edges of $G$ and 
less that $2\e |E(H)|$ edges of $H$ such that
\begin{itemize}
\item In the remaining graphs
$G'$ and $H'$, all components have size at most $K$.
\item $\sum_{S, |V(S)|\leq K} |c^{G'}_S-c^{H'}_S|<\delta,$
\end{itemize}
where $C^{G'}_S$ is the set of points that are in a component of $G'$
isomorphic to $G'$, and $c^{G'}_S=\frac{|C^{G'}_S|}{|V(G')|}$.
\end{theorem}
Thus, according to the Equipartition Theorem if a graph $H$ is
statistically close to a planar graph $H$, then $G$ can be made
planar by removing a small of amount of edges. This means exactly that
the planarity property is testable among bounded degree graphs (see \cite{BSS}).
The analogue of Theorem \ref{equi} was proved in \cite{ElekRSA} for
graph classes of subexponential growth. Using Theorem \ref{equi}, we will 
prove that if a hyperfinite graph sequence
converges then it converges locally-globally.

The following consequence of the Equipartition Theorem was proved
by Newman and Sohler \cite{NS}
 (based on the work of Hassidim, Kelner, Nguyen and Onak \cite{HKNO} ))
This result can be viewed as the finitary version of 
Theorem \ref{harmadiktetel}. 
\begin{theorem}\label{newman}
Let $\cP\subset\Gd$ be a hyperfinite family. 
Then for any $\delta>0$, there exists 
$f(\delta)>0$ such that if for a graph $G\in \cP$ and $H\in \Gd$, $|G|=|H|$,
$d_{stat}(G,H)<f(\delta)$ then $G$ and $H$ are $\delta$-close, that is
we have a bijection $\rho:V(G)\to V(H)$ such that
$$|\rho^{-1}E(H)\triangle E(G)|<\delta n\,.$$
\end{theorem}
It immediately follows from Theorem \ref{newman}, that graph isomorphism is
testable for hyperfinite graph families. 
Consequently, every reasonable property and parameter are testable for 
hyperfinite graph families (see \ref{newman} for definitions of testability).
Similar testability results were proved in \cite{ElekRSA} in case
of graph families of subexponential growth.

\vskip 0.2in
\noindent
{\bf Acknowledgement:} The author thanks the Mittag-Leffler Institute, where
parts of the paper were written, for their hospitality.

\section{Kaimanovich's Theorem revisited}
The goal of this section is to generalize a result of Kaimanovich \cite{Kai}.
First we prove a statement that is missing from \cite{Kai}, but seems to be
implicitly accepted in the paper.
\begin{definition} A graphing $\cG$ has Property $A$ if for every induced Borel
  subgraphing $\T\subseteq \cG$, almost every component have zero isoperimetric
  constant.
\end{definition}
\begin{definition} A graphing $\cG$ has Property $B$ if the following condition
is satisfied. 
For any $\e>0$, every induced subgraphing $\T$ contains a subgraphing $\SSS$ 
that
intersects almost every components of $\T$ and all the components of 
$\SSS$ have
isoperimetric constant less than $\e$ in $\T$.
\end{definition}
Note that if $F\subset \T$ is a finite subgraph, then its isoperimetric
constant
is defined as
$$i(F):=\frac{|\partial_E F|}{|F|}\,,$$
where $\partial_E F$ is the set of edges $e=(x,y)$ such that $x\in F$ and
$y\notin F$. The isoperimetric constant of an infinite graph is the infimum
of the isoperimetric constants of its finite subgraphs. An induced subgraphing
$\T$ of $\cG$ is a Borel graphing on a Borel subset $Y$ of $X$ such that
if $p,q\in Y$ is adjacent in $\cG$ then they are adjacent in $\T$ as well.
\begin{proposition}\label{liba}
For a graphing $\cG$ of vertex degree bound $d$ the two properties above are
equivalent.
\end{proposition} 
\proof
We only need to prove that Property $A$ implies Property $B$. Let $T\subseteq
\cG$ be a subgraphing satisfying the condition of Property $A$. We construct 
$\SSS\subset \T$ inductively. Let $\SSS_{n-1}\subset\T$ be the subgraphing
constructed after the $n-1$-th step consisting of finite components having
isoperimetric constants less than $\e$. Now let us consider a Borel coloring
$$\phi_n:X\to C_n=\{a_1,a_2,\dots, a_{q_n}\}$$
by finitely many colors such that $\phi_n(x)\neq \phi_n(y)$ if $d_{\cG}(x,y)\leq
2n+2$. Such coloring exists by \cite{KST}. Let $A_1=\phi^{-1}_n(a_1)$ be the
first color-class. For $x\in A_1$ let $K_x^1$ be the set of finite subsets $F$
in $B_n(x)$ containing $x$, having isoperimetric constant less than $\e$ and
such that $F\cap B_2(S_{n-1})=\emptyset$. Note that $B_2(L)$ is the
$2$-neighborhood of the set $L$.

We use the standard ordering trick and
suppose that $X=[0,1]$. Let us order $K_x^1$ the following way.
\begin{itemize}
\item If $|A|<|B|$, then $A<B$.
\item If $|A|=|B|$, then $A<B$ provided that
$$\min_{a\in A} a < \min_{b\in B}b \,.$$
\end{itemize}
Let $R_x^1$ be the smallest element of $K^1_x$. Then $\cup_{x\in A_1} R^1_x$
is a Borel set. Now let $A_2=\phi^{-1}_2(a_2)$ be the second color-class.
For $x\in A_2$ let $K_x^2$ be the set of finite subsets $F$
in $B_n(x)$ containing $x$, having isoperimetric constant less than $\e$ and
such that $F\cap B_2(S_{n-1}\cup\bigcup_{x\in A_1} R^1_x)=\emptyset$.
Again, we consider the smallest element in $K_x^2$. Then
 $\cup_{x\in A_2} R^2_x$ is a Borel set. Inductively, we define
the Borel sets $\cup_{x\in A_i} R^i_x$ and finally we define
$$S_n=S_{n-1}\cup(\bigcup_{i=1}^{q_n}\bigcup_{x\in A_i} R^i_x)\,.$$
Then $S_n$ also consists of components having isoperimetric constant less than
$\e$. Now we prove that $\SSS=\cup^\infty_{n=1}S_n$ intersects almost all
components of $\T$. Let $Z\subset \T$ be a component of isoperimetric constant
zero and let $F\subset Z$ be a finite subset of isoperimetric constant less
than $\e$. Let $F\subset B_n(x)$ for some $x\in F$. Then the only reason for
not to choose $F$ as some $R^i_x$ in the $n$-th step is that we choose another
subset
$G\subset  Z$ with isoperimetric constant less than $\e$. This shows that
Property $A$ implies Property $B$. \qed
\vskip 0.2in
\noindent
Now we are ready to state and prove Kaimanovich's Theorem.
\begin{proposition}[Kaimanovich's Theorem]
For a graphing $\cG$ of vertex degree bound $d$ the following two statements
are equivalent.
\begin{enumerate}
\item $\cG$ is hyperfinite.
\item For any subgraphing $\T\subseteq\cG$ of positive measure almost all the
  components have isoperimetric constant zero.
\end{enumerate}
\end{proposition}
\proof First we show that (2) implies (1). Let us suppose that $\cG$ satisfies
the second condition. Let $\cG=\T_0$ and $\SSS_0$ be a Borel subset of positive
measure consisting of finite components with isoperimetric constant less than
$\e>0$. Such set exists by Proposition \ref{liba}. Let $E_0$ be the set of 
edges pointing out of $\SSS_0$. Then
$\mu_E(E_0)\leq \e\mu(\SSS_0)\,.$
Remove $E_0$ from $\T_0$ along with the subgraphing $\SSS_0$. Let us denote the
resulting subgraphing by $\T_1$. Note that
$\mu(\T_1)<\mu(\T_0)$, where $\mu(\T_1)$ denote the measure
of the vertex set of $\T_1$. Now we proceed by transfinite induction.
Suppose that $\T_\alpha$ is constructed for some countable ordinal and 
$\mu(\T_\alpha)>0$. Let $\SSS_\alpha$ be a Borel set of positive measure
consisting of finite components with isoperimetric constant less than
$\e>0$ in $\T_\alpha$. Again, let $E_\alpha$ be the set of 
edges pointing out of $\SSS_\alpha$. Then
$\mu_E(E_\alpha)\leq \e\mu(\SSS_\alpha)\,.$
Remove $E_\alpha$ from $\T_\alpha$ along with the subgraphing $\SSS_\alpha$.
 Let us denote the
resulting subgraphing by $\T_{\alpha+1}$. Then $\mu(\T_{\alpha+1})<\mu(\T_{\alpha})$.
For a limit cardinal $\alpha'$, let $\T_\alpha'$ be $\cap_{\alpha<\alpha'} 
\T_\alpha$.
 Since
$\mu(\T_\alpha)>\mu(\T_{\alpha+1})$, there exists a countable ordinal $\beta$
for which $\mu(\T_\beta)=0$. Let $\SSS=\cup_{\beta<\alpha} \SSS_\alpha$ and 
$M=\cup_{\beta<\alpha} E_\alpha$. Clearly, $\mu_E(M)<\epsilon$. Hence, by
removing $M$ and $T_\alpha$ from $\cG$ we obtain a graphing consisting of finite
components. This implies the hyperfiniteness of $\cG$. 
\vskip 0.2in
\noindent
Now let us prove that (1) implies (2). Suppose that
$\cG$ has a subgraphing $\T$ of positive measure such that
the measure of points $p$ for which the component $Z_p$ has positive
isoperimetric constant is not zero. Then there exists $\delta>0$ and a Borel
subgraphing $\T_\delta\subset \T$ of positive measure such that all the
components of $T_\delta$ have isoperimetric constants at least $\delta$.
Now suppose that $\cG$ is hyperfinite. Let $F$ be a Borel set of edges such that
 $\mu_E(F)<
\frac{\delta |\mu(T_\delta)|}{10}$ and $\SSS=\cG\backslash F$ consists
of finite components. Let $K$ be a component of $\SSS$. Then by our condition,
there exist at least $\delta|T_\delta\cap K|$ edges pointing out of $K$.
This gives us an estimate for the edge density of $F$
$$\mu_E(F)> \delta |T_\delta|\,,$$ leading to a contradiction.
 \qed
\vskip 0.2in
\noindent
Kaimanovich's Theorem will be applied in our paper using
the following corollary.
Let $\cG(X,\mu)$ and $\cH(Y,\nu)$ be graphings. A surjective map 
$\pi:X\to Y$ is a factor map (that is $\cH$ is a factor of $\cG$) if
\begin{itemize}
\item $\pi$ is measure preserving, that is for any Borel set $A\subseteq Y$,
$\mu(\pi^{-1}(A))=\nu(A)$.
\item For almost all $x\in X$, $\pi$ is a graph isomorphism restricted on the 
component of $x$.
\end{itemize}
\begin{proposition} \label{factor}
If $\cH$ is a factor of $\cG$, then $\cH$ is hyperfinite if and only if
$\cG$ is hyperfinite.
\end{proposition}
\proof
First suppose that $\cH$ is hyperfinite and
$W$ is a Borel set of the edges of $\cH$ such that
$\nu_E(W)<\epsilon $ and all the components of $E(\cH)\backslash W$
have size at most $K$. Then $\mu_E(\pi^{-1}(W))<\epsilon$ and all the 
components of
$E(\cG)\backslash \pi^{-1}(W)$ have size at most $K$. Hence $\cG$ is 
hyperfinite. 

\vskip 0.2in
\noindent
For the converse statement, suppose that $\cH$ is not hyperfinite.
Then by Kaimanovich's Theorem, there exists a subgraphing $\cT\subseteq \cH$
such that {\it not} almost all its components have zero isoperimetric
constant. Then $\pi^{-1}(\cT)$ is a subgraphing of $\cG$ witnessing
the non-hyperfiniteness of $\cG$. \qed

\vskip 0.2in
{\bf Note:} Let $\Gamma$ be a finitely generated amenable group acting
freely on the standard Borel space $(X,\mu)$ preserving the probability
measure. Then the graphing of the action is hyperfinite.
The standard proof of this fact is given by the Ornstein-Weiss quasi-tiling
construction \cite{OW}. However, a very short proof can be obtained by
Kaimanovich's Theorem. Without claiming any originality, we provide a proof
for completeness.

\vskip 0.2in
\proof 
Let $S$ be a symmetric generating system and $\cG(X,\mu)$ the graphing
of the action. Suppose that $\cG$ is not hyperfinite.
Then it contains a subgraphing $\cT$, $V(\tau)>0$, such that the isoperimetric
constants of
all the components of  $\cT$ are larger than a certain positive constant
$\delta$.
Indeed, if $\cT$ is a subgraphing with components of positive isoperimetric
constants and $\cT_\delta$ is the subgraphing consisting of components
having isoperimetric constant larger than $\delta$, then
$\cup_{\delta} \cT_\delta=\cT$.
Let $\{F_n\}^\infty_{n=1}$ be a F\o lner sequence in $\Gamma$.
By the invariance of the measure,
$$\int_X |F_nx\cap V(\cT)| d\mu= |F_n| \mu(V(\cT))\,.$$
Hence, we have a sequence of points $\{x_n\}^\infty_{n=1}\subset X$
such that
$$\frac{|F_nx\cap V(\cT)|}{|F_n|}>\mu(V(\cT))>0\,.$$
Therefore the isoperimetric constants of the induced subgraphs
$\{[F_nx\cap V(\cT)]\}^\infty_{n=1}$ tend to $0$, leading to a contradiction.
\qed

\section{Canonical limits} \label{canonical}
This section is of rather technical nature.
\subsection{The Benjamini-Schramm limit measure}
First let us recall the notion of the 
Ben\-jamini-Schramm limit measure 
construction. Let $\grd$ be the set of all connected, rooted, countable 
graphs up to rooted graph isomorphisms.
One can introduce a metric on $\grd$ by setting
$$d(X,Y)=2^{-r}\,,$$
where $r$ is the largest integer such that the $r$-balls around the roots
of $X$ resp. $Y$ are isomorphic. The metric space $\grd$ is compact.
Note that for all $r\geq 1$ and $\alpha\in\urd$, $T(\alpha)\subseteq \grd$,
that is the set of all graphs such that the $r$-ball around their roots is 
isometric to $\alpha$ is a clopen set.
Now let $\hg=\{G_n\}\subset \Gd$ be a convergent graph sequence. Then
$$\mu_{\hg}(T(\alpha))=\lim_{n\to\infty} p(G_n,\alpha)$$
defines a Borel probability measure on $\grd$.
This measure is called the Benjamini-Schramm limit measure 
(a so-called unimodular measure, see \cite{AL})
We say that $X,Y\in\grd$ are adjacent if there is a neighbouring vertex $y$ of
the root of $X$ such that $Y$ is rooted isomorphic to the underlying graph
of $X$ with root $y$. In this way, $\grd$ is equipped with a Borel graph
structure. However, the following example shows that $(\grd,\mu_{\hg})$
is not necessarily a graphing.
\begin{example}
Let $G_n$ be the graph obtained from the line graph $L_n$ of length $n$
by adding two leaves for each vertex. Then $\hg=\{G_n\}^\infty_{n=1}$ is
a convergent graph sequence. The limit measure is concentrated to two
points $a$ and $b$ such that $\mu_{\hg}(a)=1/3$ and $\mu_{\hg}(b)=2/3$. Hence
$(\grd,\mu_{\hg})$ is not a graphing.\end{example}

\subsection{$B$-graphs} It was observed by Aldous and Lyons (Example 9.9 
\cite{AL}) that for each 
unimodular measure, one can construct a marked network, which is a graphing.
This should be thought as the Bernoulli space 
of the unimodular measure. So let us recall the notion of $B$-graphs 
from \cite{ElekRSA}. 
This is an explicite realization of the Aldous-Lyons marked network 
construction.
Let $B$ be the set $\{0.1\}^{\NN}$ with the standard product measure. 
Then $\gdb$ is the set of all finite
simple graphs of vertex degree bound $d$  with vertices colored by $B$ 
(up to colored isomorphisms). 
These objects are called $B$-graphs.
Let $\urdb$ be the set of all rooted $r$-balls with vertices colored
with $\{0,1\}$-strings of length $r$. If $G\in \gdb$, $\beta\in \urdb$ and
$x\in V(G)$ then $x\in T(G,\beta)$ if the rooted $r$-ball around $x$ is
isomorphic to $\beta$, when one restricts the color of the vertices to
the first $r$-digits. Set $p(G,\beta):=\frac{T(G,\beta)}{|V(G)|}$.
Again, we can define the convergence of $B$-graphs.
The sequence of $B$-graphs $\{G_n\}^\infty_{n=1}$ is convergent if
for any $r\geq 1$ and $\beta\in\urdb$, $\limn p(G_n,\beta)$ exits.

\noindent
The corresponding limit objects are measures on $\grdb$, the space of
connected, rooted, countable, $B$-colored graphs. The reason we introduced 
the notion of $B$-graphs is that using them one can construct canonical limit 
graphings of standard (colorless) convergent graph sequences.
Let us recall the construction from \cite{ElekRSA}.
Let $\hg=\{G_n\}^\infty_{n=1}\subset \Gd$ be a convergent 
graph sequence. Let us color
the vertices of the graphs in the sequence randomly, independently, by elements
of the probability measure space $B$. Then with probability $1$, the resulting
$B$-colored graph sequence will be convergent to the same measure 
$\mu^B_{\hg}$ on $\grdb$. This measure is the Bernoullization of the
Benjamini-Schramm limit measure of the original graph sequence in the sense of
Aldous and Lyons. 
 Then
\begin{itemize}
\item Then $\cG(\grdb,\mu^B_{\hg})$ is a graphing that we denote by
$\cG_{\hg}$.
\item For $\mu^B_{\hg}$-almost all $x\in\grdb$ the orbit of $x$ in
$\cG_{\hg}$ is isomorphic to $x$ as rooted $B$-graphs.
\end{itemize}
We call $\cG_{\hg}$ the canonical limit graphing of $\hg$.
\subsection{Edge-colored graphs and $B$-graphs}
Finally, we need a little bit more complicated construction of the same 
genre as the ones above. Let $C\Gd$ be the set of simple graphs of vertex
degree bound $d$ with proper edge-colorings by ${d+1\choose 2}$ colors.
Recall that a coloring is proper if incident edges are colored 
differently. Similarly, we can consider $C\gdb$ the set of $B$-graphs
of vertex degree bound $d$ with proper edge colorings by 
${d+1 \choose 2}$ colors.
Again, we can define the convergence of edge-colored graphs resp. edge-colored 
$B$-graphs together with compact metric spaces $C\grd$ resp. $C\grdb$, the
spaces of properly edge-colored rooted graphs resp. properly edge-colored
rooted $B$-graphs (by ${d+1\choose 2}$ colors).
Also, the limits of convergent sequences are the appropriate measures
on $C\grd$ resp. $C\grdb$. One should note that there exists a natural
$\Gamma$-action on graphs properly edge-colored by ${d+1\choose 2}$ colors,
where $\Gamma$ is the ${d+1\choose 2}$-fold free product of cyclic
groups of order $2$. Also, $\Gamma$ acts on $C\grd$ resp. on $C\grdb$
by homeomorphisms. Let $\hat{H}=\{H_n\}^\infty_{n=1}\subset C\Gd$ be 
a convergent graph sequence and $\mu_{\hat{H}}$ be the limit measure
on $C\grd$. Then the Borel probability measure $\mu_{\hat{H}}$ is 
invariant under the natural $\Gamma$-action. Similarly to the colorless
case one can construct the canonical limit measure $\mu^B_{\hat{H}}$ on
$C\grdb$ as well.

\section{The Oracle Method}
The essence of the oracle method is that it enables us to 
construct subsets of finite graphs using one single subset of
$\grdb$.
The Oracle Method is strongly related to the notion of
randomized distributed algorithms. Suppose that  a subset
$A\subseteq U^{r,B}_d$ is given. 
Say, we have a finite graph $G$
of degree bound $d$. We color the vertices of $G$ random uniformly
with $\{0,1\}$-strings of length $r$.
Then we construct a subset $V_A\subseteq V(G)$ the following way.
If the $r$-ball around $v\in V(G)$ is colored-isomorphic to
an element of $A$, let $v\in V_A$. Otherwise, $v\notin V_A$.
The only reason we need colorings is that
we can use the colors to ``break ties'' in the case of symmetries.
If $G$ is a transitive graph, distributed algorithms without
randomization can produce only the empty set and $V(G)$ itself.

Now let $x\in\grdb$ and $\mu^B_{\hg}$ be a limit measure. As
it pointed out in Section \ref{canonical}, the measure
$\mu^B_{\hg}$ is concentrated on countable graphs with ``broken'' symmetries
that is on graphs for which all the vertex colors are different.
In this case, the component of $x$ in the Borel graphing $\cG_{\hg}$
is isomorphic to the underlying graph of $x$.
Of course, if the underlying graph of $x$ is transitive and all the
vertex colors are identical, then the component of $x$ is just one single 
vertex. In this case, we lose all the information about the graph structure
of $x$. If the colors on the $r$-ball around the root of $x$ are different,
then we know at least that the $r$-ball around the root of $x$ and the
$r$-ball around $x$ in the graphing are isomorphic.

In order to handle the color issue, we need a simple
variation of $U^{r,B}_d$. Let $s>r$ be an integer. Then
$U^{r,s,B}_d$ is the set of $r$-balls with vertices
colored by  $\{0,1\}$-strings of length $s$. Obviously, $U_d^{r,r,B}=
U^{r,B}_d$. Let $W^{r,s,B}_d\subset U^{r,s,B}_d$ be the set of balls
for which the vertex colors are all different.
Let $V^{r,s,B}_d= U^{r,s,B}_d\backslash W^{r,s,B}_d$.
The following lemma is an easy consequence of the law of large numbers
and is left to the reader.
\begin{lemma}
For any $\delta>0$ and $r\geq 1$ there exists $s>r$ such that
$$\mu_{\hg}\left( \cup_{\alpha\in V_d^{r,s,B}} \mu^B_{\hg} (T(\alpha)\right)
<\delta $$
for any convergent graph sequence $\hg$.
\end{lemma}
\begin{proposition} \label{egyik}
Let $\hg=\{G_n\}^\infty_{n=1}\subset G_d$ be 
a convergent graph sequence and let $\cG_{\hg}$ be its
canonical limit graphing. If
$\cG_{\hg}$ is hyperfinite, then $\hg$ is always a hyperfinite
family.
\end{proposition}
\proof Fix a constant $\delta>0$.
Let $N\subset \grdb$ be a Borel subset such that $\mu_{\hg}(N)<\delta$
and if we remove the edges incident to the vertices
in $N$, then all the components of the resulting subgraphing have size at
most $K$.
The goal is to prove that the graphs inherit this property. That is, 
if $n$ is large 
enough, then there exists $P_n\subset V(G_n)$, $|P_n|<\delta|(V(G_n)|$ such 
that if we remove the edges of $G_n$ incident
to the vertices of $P_n$, the resulting graph $G'_n$ has components
of size at most $K$. The following approximation lemma is
the key of the proof of Proposition \ref{egyik}.
\begin{lemma}\label{appro}
Let $\hg$ and $K$ be as above. Then there exist
integers $s>r>K$ and a subset $V^{r,s,B}_d\subset A\subset U^{r,s,B}_d$ with
the following property.
\begin{itemize}
\item $\mu_{\hg}(N_A)<\delta$, where $N_A=\cup_{\beta\in A} T(\beta)$ \,.
\item If we remove the edges of $\cG_{\hg}$ incident to points in $N_A$, 
the components
of the resulting subgraphing $\cT$ are of size at most $K$.\end{itemize}
\end{lemma} One can interpret the lemma in the following way. The 
hyperfiniteness of $\cG_{\hg}$ can be witnessed by the removal of edges
incident to ``nice'' subsets.
First, let us show that the lemma implies Proposition \ref{egyik}.
Let $t>s$, $t>2r$ be an integer such that
$$\mu^B_{\hg}(\cup_{\beta\in V^{2r,t,B}_d }T(\beta))< \delta-\mu_{\hg}(N_A)\,.$$
Take a random coloring of the vertices of the graphs $\{G_n\}^\infty_{n=1}$
by $B$. Let $H_n\subset V(G_n)$ be the set of vertices $x$
such that either $B_r(x)\in A$ or
$B_{2r}(x)\in V_d^{2r,t,B}\,.$
Remove the edges incident to $H_n$. Then in
the resulting graph $G'_n$ the maximal component size is at most $K$.
Indeed, suppose that there is a component of size greater than $K$
and $v\in K$. Then $B_{2r}(v)\in W^{2r,t,B}_d$, thus the
$2r$-ball around $v$ in $G_n$ is isomorphic to the $2r$-ball round
the point $z\in \grdb$, where $z\in T(\alpha)$, $\alpha\sim B_{2r}(v)\,.$
Observe, that by our construction, $B_r(x)\cap G_n'$ must be a subgraph of 
$B_r(z)\cap\cT$. Since the later graph does not contain components
of size larger than $K$, neither does $B_r(x)\cap G_n'$.
Therefore the maximal component size in $G_n'$ is at most $K$.
Now Proposition \ref{egyik} follows from the fact that for any $\alpha\in 
U^{r,s,B}_d$, $\limn p(G_n,\alpha)=\mu_{\hg}(T(\alpha))$ with
probability one. \qed

\vskip 0.2\in
Now let us prove Lemma \ref{appro}.
Let $H\subset \grdb$ be a Borel subset, $\mu_{\hg}(H)<\delta$
such that if we remove the edges incident to $H$, the remaining
components have size at most $K$. Since sets in the form $N_A$, where
$A\subset U^{l,B}_d$ for some $l>K$
generate the Borel sets of $\grdb$
we have a sequence $\{N_{A_l}\}^\infty_{l>K}$ such
that
\begin{equation} \label{e61}
\lim_{l\to\infty} \mu^B_{\hg}(N_{A_l}\triangle H)=0\,.
\end{equation}
Let $\cT_l$ be the subgraphing obtained from the Borel graphing $\cG_{\hg}$
by removing the edges incident to $N_l$. Let $X_l$ be the set of points
in $\grdb$ that are in a component of $T_l$ larger than $K$.
Observe that $\lim_{l\to\infty}\mu^B_{\hg}(X_l)=0\,.$
Pick $s(l)>2l$ in such a way that $\lim_{l\to\infty}\mu_{\hg}(P_l)=0$,
where $P_l=\cup_{\alpha\in V^{2l,s(l),B}_d} T(\alpha)\,.$
Let $Q_l=\cup_{\beta} T(\beta)$, where
the index $\beta$ runs through all elements of $W^{2l,s,B}_d$ 
such that the root of $\beta$ is contained in a component of $\cT_l$
larger than $K$. Note that it is meaningful, since by looking at
the $2l$-neighborhood of a vertex we can decide whether it is
contained in a component of $\cT_l$ larger than $K$.
Since $Q_l\subseteq X_l$, $\lim_{l\to\infty}\mu^B_{\hg}(Q_l)=0\,.$
Hence if $l$ is large enough then
$N_A=N_{A_l}\cup P_l\cup Q_l$ satisfies the conditions of the lemma 
(with $r=2l)\,.$ \qed
\vskip 0.2in
\noindent
Now we prove the converse
of Proposition \ref{egyik}.

\begin{proposition} \label{masik}
Let $\hg=\{G_n\}^\infty_{n=1}$ be a hyperfinite convergent graph sequence.
Then the canonical limit $\cG_{\hg}$ is hyperfinite.
\end{proposition} 
\proof
As in the previous sections, let us color the vertices in the
graph sequence randomly by $B$.
Now we construct a second $B$-coloring of the vertices.
The $k$-th digit of the second $B$-color of $x\in V(G_n)$
is given the following way.
Let $C_k$ be an integer such that
for any $n\geq 1$ there exists a subset $H_{n,k}\subset V(G_n)$, with
${|H_{n,k}|}{|V(G_n)|}<1/k$
such that if we remove the edges incident to $H_{n,k}$, the
components in the remaining graph $G_{n,k}$ have size at most $C_k$.
Let $0$ be the $k$-th digit of $x$ if $x\in H_{n,k}$, otherwise
let the $k$-th digit be $1$. This way we constructed a coloring
of the graphs by $B^2$. Note that for convergent $B^2$ colorings we have
limit measures on $\grdbk$ completely analogously to $B$-colorings.
We cannot say that the $B^2$-colored graphs constructed above are convergent 
(as colored graphs). However, we have a convergent subsequence by compactness.
Let ${\mu}^{B^2}_{\hg}$ be the associated limit measure on
$\grdbk$. 
Then, $\pi:\grd^{B^2}\to\grdb$ is a factor map, where $\pi$ forgets
the second coordinate.
Now let us observe that the graphing $\cG(\grd^{B^2}, {\mu}^{B^2}_{\hg})$
is hyperfinite. Indeed, the Borel set of vertices with $0$ as the $k$-th
digit of their second $B$-coordinate has ${\mu}^{B^2}_{\hg}$
measure less than $1/k$. Also, if we remove the edges 
incident to this set the remaining graphing have components of size
at most $C_k$. By Proposition \ref{factor}, our proposition follows. \qed

\section{The proof of Theorem \ref{elsotetel} and Theorem \ref{masodiktetel}}
We will show slightly more. Let $\cH(X,\mu)$ be an arbitrary graphing 
with vertex degree bound $d$,
We can consider the associated unimodular measure the following way 
(\cite {AL},\cite{Artemenko}). 
For each point $x\in X$ let $\pi(x)\in \grd$ be the
component of $x$ in $\cH$ with $x$ as the root. Then the measure $\pi_*(\mu):=
\mu_{\cH}$
is unimodular (see also Corollary 6.10 \cite{Artemenko}).
We can consider the Bernoulli measure $[\pi^*(\mu)]_B:=\mu^B_{\cH}$ 
on $\grdb$ (see Section
\ref{canonical}) and the corresponding graphing $\cG(\grdb,\mu^B_{\cH})$.
If $\cG$ is weakly equivalent to $\cH$, then the associated
Bernoulli measures and the corresponding graphings are the same. Hence by
Proposition \ref{factor}, the following lemma immediately implies 
Theorem \ref{masodiktetel}.
\begin{lemma}
There exists a graphing $\cK$ such that $\cH$ and $\cG(\grdb,\mu^B_{\cH})$
are both factors of $\cK$.
\end{lemma}
\proof
First let us note that if the measure $\mu_{\cH}$ is concentrated
on rooted graphs without rooted automorphisms, then 
$\cH\to \cG(\grd,\mu_{\cH})$ is already a factor map. In this case, the proof
of Theorem 2. would end here. We use the Bernoullization only to
handle the symmetries.
This is the point in our paper, where we use the edge-colorings.
By a result of Kechris, Solecki and Todorcevic \cite{KST}
one can color the vertices of a Borel graphing of vertex
degree bound $d$ properly with $d+1$-colors in a Borel way.
This vertex coloring gives us a Borel edge coloring of $\cH$
with $d+1 \choose 2$-colors. The color of an edge between a
vertex colored by $a$ and a vertex colored by $b$ will be colored
by $(a,b)$.
As it was mentioned in Section \ref{canonical},
the coloring defines a Borel $\Gamma$-action on $X$, where
$\Gamma$ is the free product
of $d+1 \choose 2$ cyclic groups of order $2$.
Again, we have the natural $\Gamma$-equivariant map
$\pi_C:X\to C\grd$.
We denote $(\pi_C)_{\star}(\mu)$ by $\mu^C_{\cH}$.
Now let us consider the Bernoullization of $\mu^C_{\cH}$ on
$C\grdb$, $\mu^{C,B}_{\cH}$.  We have two factor maps
\begin{itemize} \item
$\rho:\cG(C\grdb, \mu^{C,B}_{\cH})\to \cG(\grdb,\mu^B_{\cH})$
the map forgetting the edge-colorings.
\item
$\zeta: \cG(C\grdb, \mu^{C,B}_{\cH})\to \cG(C\grd, \mu^C_{\cH})$
the map forgetting the vertex-colorings.
\end{itemize}
Note that $\pi_C$ and $\zeta$ are both $\Gamma$-equivariant maps
so, we can consider their relative independent joining \cite{Gla}
over $C\grd$. This gives us a new $\Gamma$-action on a space $Y$, with
graphing $\cK$. By the joining construction,
both $\cH$ and $\cG(C\grdb, \mu^{C,B}_{\cH})$ are factors of $\cK$.
On the other hand, $\cG(\grdb,\mu^B_{\cH})$ is a factor of
$\cG(C\grdb, \mu^{C,B}_{\cH})$. Thus both $\cH$ and
$\cG(\grdb,\mu^B_{\cH})$ are factors of $\cK$. Hence the lemma follows.
\qed

\vskip 0.2in
Now let us observe that Theorem \ref{elsotetel} immediately follows from
Theorem \ref{masodiktetel} by Proposition \ref{egyik} and Proposition
\ref{masik}.

\section{Equipartitions}
\subsection{The Transfer Theorem}
The Transfer Theorem is one of the basic applications of the
Oracle Method. 
\begin{theorem}[Transfer Theorem]
Let $\hg=\{G_n\}^\infty_{n=1}\subset \Gd$  be a convergent graph sequence.
Let $\cH\subseteq \cG_{\hg}$ be a subgraphing 
(note that it means that $V(\cH)=\grdb$).
Then there exist subgraphs $H_n\subseteq G_n$, $V(G_n)=V(H_n)$ such that
$\{H_n\}^\infty_{n=1}$ converges to $\cH$.
\end{theorem}
\proof
Recall that $\cG_{\hg}=\cG(\grdb,\mu_{\hg})$. Also, for $\mu_{\hg}$-almost
all elements $x$ of $\grdb$ the vertices of $x$ are $B$-colored differently.
Let us call such vertex $x$ typical. Thus the orbit of a typical
vertex is isomorphic to the graph represented by $x$ in $\grdb$.
How can we encode the edge set of $\cH$ ?
A symbol $\sigma$ consists of the following data.
A number $0\leq k \leq d$, the degree of the symbol, and a
subset $\{a_1< a_2<\dots <a_l\}$ of $\{1,2,\dots k\}$, where
$l\leq k$. For any edge $(x,y)\in E(\cG)$, for which $x$ is typical we
have an ``edge code'' which is $s$ is $y$ is the $s$-th neighbour of $x$
with respect to the lexicographical ordering of $B$.
If $x$ is a typical vertex, then its position in $\cH$ can be described
by the the symbol $\sigma=(k,a_1,a_2,\dots,a_l)$, where $k$ is the degree
of $x$ in $\cG$ and $a_i$ is the edge code of the $i$-th neighbor of $x$
in $\cH$ in the lexicographical ordering of the $B$-colors. We denote by
$\cH_\sigma$ the Borel set of typical vertices $x$ with $\cH$-position
symbol $\sigma$. Let $E(\cH_\sigma)$ be the set of edges in $\cG$ incident
to an element of $\cH_\sigma$. Then, $E(\cH)=\cup_\sigma E(\cH_\sigma)$.
Note that the sets $\cH_\sigma$ are disjoint. 

\vskip 0.2in
\noindent
As in Lemma \ref{appro} let $A^l_\sigma\subset U^{l,B}_d$ such that
\begin{itemize}
\item The degree of $z\in A^l_\sigma$ is the degree of $\sigma$.
\item $\lim_{l\to\infty} \mu^B_{\hg}(N_{A^l_\sigma}\triangle \cH_\sigma)=0\,.$
\end{itemize}
We also suppose that the sets $A^l_\sigma$ are disjoint.
Then one can consider the
approximating graphings $\cH^l$ 
$$E(\cH^l)=\cup_{\sigma} E([A^l_\sigma]_\sigma)\,,$$
where $E([A^l_\sigma]_\sigma)$ is the set of edges $(z,w)$ such that
$z\in A^l_\sigma$ and the ``edge code'' of $w$ belongs to $\sigma$.
Then $\lim_{l\to\infty} \cL(\cH^l)=\cL(\cH)$.
Therefore it is enough to prove the Transfer Theorem for the subgraphings
$\cH^l$. We construct the subgraphs $\{H_n\}^\infty_{n=1}$ the following way.
First, we B-color the vertices of the graphs $G_n$ randomly to obtain
the graph $G_n^B$.
Then for each vertex $v\in G_n$ we check the $l$-neighborhood of $v$.
If for some $\sigma$, $B^l(v)\in A^l_{\sigma}$ then using the symbol $\sigma$
and the $B$-coloring we choose the appropriate edges of $G_n$ incident to
$v$. In this way we obtain the subgraph $H_n$. 
The following lemma finishes the proof of our theorem.
\begin{lemma}
$\{H_n\}$ converges to $\cH_l$ with probability $1$.
\end{lemma}
\proof Let $r>0$ and $\beta\in U^{r,B}_d$. It is enough to see
that
\begin{equation} \label{transeq}
\lim_{n\to\infty} p(H_n,\beta)=\mu^B_{\hg}(T(\cH^l,\beta))
\end{equation}
with probability $1$.
Let $z\in\grdb$, $z\in T(\gamma), \gamma\in W^{r+l,s,B}_d$.
Then $\gamma$ determines whether $z\in T(\cH^l,\beta)$ or not.
We denote by $W^{r+l,s,B}_{d,1}$ the set of $\gamma$'s where
$z\in T(\cH^l,\beta)$. So, if $v\in T(G_n^B,\gamma)$, $\gamma\in W^{r+l,s,B}_{d,1}$
then $v\in T(H_n,\beta)$ and if $v\in T(G_n^B,\gamma')$,
$\gamma'\in W^{r+l,s,B}_d\backslash  W^{r+l,s,B}_{d,1}$, then
$v\notin T(H_n,\beta)$. Hence we have the following estimates.
$$\sum_{\alpha\in W^{r+l,s,B}_{d,1}} p(G^B_n,\alpha)\leq
p(H_n,\beta)\leq \sum_{\alpha\in W^{r+l,s,B}_{d,1}} p(G^B_n,\alpha) +
\sum_{\alpha'\in V^{r+l,s,B}_d} p(G_n,\alpha)\,.$$
and
$$\sum_{\alpha\in W^{r+l,s,B}_{d,1}} \mu^B_{\hg}(T(\alpha))\leq
\mu^B_{\hg}(T(\cH^l,\beta))\leq \sum_{\alpha\in W^{r+l,s,B}_{d,1}}
\mu^B_{\hg}(T(\alpha))+
\sum_{\alpha'\in V^{r+l,s,B}_d}\mu^B_{\hg}(T(\alpha))\,.$$
Since 
$$\lim_{s\to\infty} \sum_{\alpha'\in V^{r+l,s,B}_d}\mu^B_{\hg}(T(\alpha))=0$$
the lemma follows from the fact that $\{G_n\}^\infty_{n=1}$ is a 
convergent sequence. \qed
\subsection{The Uniformicity Theorem}
Let $\cP\subset \Gd$ be a hyperfinite family.
Denote by $L\cP$ the set of graphings that are limit graphings
of sequences in $\cP$. By Theorem \ref{elsotetel}, the elements
of $L\cP$ are hyperfinite graphings. The Uniformicity Theorem states
that $L\cP$ is a uniformly hyperfinite family of graphings.
\begin{theorem}[The Uniformicity Theorem]
Let $\cP\subset \Gd$ be a hyperfinite family then for any $\zeta>0$
there exists $K>0$ such that
for each $\cG\in L\cP$ there exists a Borel set $Z\subset E(\cG)$
of edge-measure less than $\zeta$ such that
the components of $\cG\backslash Z$ are of size at most $K$.
\end{theorem}
Let $\cH(X,\mu)$ be a hyperfinite graphing such that
all of its components have size at most $K$.
For a connected graph $S$ of size at most $K$ let $c^{\cH}_S$ be the
$\mu$-measure of points in $X$ that belong to a component
isomorphic to $S$.
Let $\{H_n\}^{\infty}_{n=1}$ be a graph sequence converging to $\cH$ and 
$C^{H_n}_S$ be
the set of vertices in $V(H_n)$ that belong to a component isomorphic to
$S$.
\begin{lemma} \label{u1}
If $\{H_n\}^\infty_{n=1}$ and $\cH$ are as above then 
$\lim_{n\to\infty} c^{H_n}_S=c^{\cH}_S.$
\end{lemma}
\proof
Let $U^{k+1}_{d,S}$ be the set of elements of $U^{k+1}_d$ that are isomorphic
to $S$. Note that these rooted balls are already in $U^k_d$. However, if
the $k+1$-ball of a vertex is in $U^{k+1}_{d,S}$ then we know that the vertex is 
in a component isomorphic to $S$.
Clearly,
$$\sum_S\sum_{\alpha\in U^{k+1}_{d,S}} \mu_{\cH}(\alpha)=1\,,$$
where $S$ is running through the isomorphic classes of
connected graphs of size at most $K$. By convergence, for any $S$ and
any $\alpha\in U^{k}_{d,S}$
$$\lim_{n\to\infty} p(H_n,\alpha)=\mu_{\cH}(T(\alpha))\,.$$
Observe that
$c^{\cH}_S=\sum_{\alpha\in U^k_{d,S}} p(H_n,\alpha)\,.$
Hence the lemma follows. \qed

\vskip 0.2in
The proof of the next lemma is basically identical to the previous one.
\begin{lemma}\label{u2}
Let $\{\cH_n\}^\infty_{n=1}$ be a sequence of
graphings such that $\lim_{n\to\infty}\cL(\cH_n)=\cL(\cH)$, where $\cH$ is
as above. Then $\lim_{n\to\infty} c^{\cH_n}_S=c^{\cH}_S,$ for any $S$, 
$|V(S)|\leq K$.
\end{lemma}

Now let $\hg=\{G_n\}^\infty_{n=1}$ be a convergent
graph sequence and let $Z\subset E(\cG_{\hg})$ be a Borel set
of edges with edge-measure less than $\e>0$,
such that the subgraphing $\cH=\cG_{\hg}\backslash Z$
consists of components of size at most $K$.
For the rest of this section we consider this subgraphing $\cH$.
We will show that if a hyperfinite graphing $\cG'$ is statistically
close to $\cG$ then it contains a subgraphing $\cH'$ of
components of size at most $K$, such that
$c^{\cH}_S$ is close to $c^{\cH'}_S$ for any connected graph $S$, $|V(S)|\leq K$.
First we formulate this statement
for finite graphs.
\begin{lemma} \label{u3}
Let $\hg,\cH,K$ be as above. Then for any $\delta>0$ there exists
$f(\delta)>0$ such that if for a 
finite graph $G$ $d_{stat}(G,\cG_{\hg})<f(\delta)$ then
$G$ contains a subgraph $H\subset G$ with components of size at
most $K$ such that
\begin{equation} \label{ronda}
|c^H_S-c^{\cH}_S|<\delta  \,\,\,\mbox{for all}\,\,\, S, |V(S)|\leq K.
\end{equation}
\end{lemma}
\proof
Suppose that the lemma does not hold. Then we have a $\delta>0$
and a sequence of finite graphs $\hat{Q}=\{Q_n\}^\infty_{n=1}$
converging to $\cG_{\hg}$ without subgraphs $H_n$ satisfying
(\ref{ronda}).
Observe that $\cG_{\hat{Q}}=\cG_{\hg}$. Thus, by the Transfer Theorem
there exists subgraphs $H'_n\subset Q_n$ converging to $\cH$.
By Lemma \ref{u1}, $\lim_{n\to\infty} c^{H'_n}=c^{\cH}_S$, for any $S$.
So, we have subgraphs $H_n\subset H'_n$ with components of size
at most $K$ such that $\lim_{n\to\infty} c_S^{H_n}=c^{\cH}_S$ for any $S$,
leading to a contradiction. \qed

\begin{lemma} \label{becslo}
Let $\cK(X,\mu)$ be a graphing such that
all of its components are of size at most $l$. Let $\{Q_n\}^\infty_{n=1}$
be a sequence of graphs converging to $\cK$.
Let $H_n\subset Q_n$ be subgraphs with components of size at most $K$
such that for all $n\geq 1$, $|c_S^{H_n}-c^{\cH}_s|<\delta/4$, where
$\cH$ is the subgraphing as above. 
Then there exists a subgraphing $\cH'\subset \cK$ with components of
size at most $K$, such that
$|c^{\cH}_S-c^{\cH'}_S|<\delta/2$ for all connected graph $S$, $|V(S)|\leq K$.
\end{lemma}
\proof
Let $L$ be a connected graph, $|V(L)|\leq l$. Let $C^{Q_n}_{S,L}$
be the set of points in $Q_n$ that are in a component
$C$ of $Q_n$ such that $C\cap L\cong S$. Set
$c^{Q_n}_{S,L}=\frac{C^{Q_n}_{S,L}}{|V(Q_n)|}$. Then
$\sum_S c^{Q_n}_{S,L}=c^{Q_n}_L$.
Pick a subsequence $\{Q_{n_k}\}^\infty_{n=1}$ such that for all $S$, $L$,
$\lim_{n\to\infty} c^{Q_{n_k}}_{S,L}=d_{S,L}$ exists.
Then $\sum_S d_{S,L}=C^{\cK}_L$.
Let $C^{\cK}_L$ be the set of points in $X$ that are in a component of $\cK$
isomorphic to $L$. Then $\mu(C^{\cK}_L)=c^{\cK}_L$. Divide $C^{\cK}_L$ into
Borel subsets such that
\begin{itemize}
\item $\mu(C^{\cK}_{S,L})=d_{S,L}$.
\item Each component of $C^{\cK}_{S,L}$ is isomorphic to $L$.
\end{itemize}
Let $\cH^{\cK}_{S,L}$ be a Borel graph on $C^{\cK}_{S,L}$, such that its
edges are edges of $\cK$ and all the components are isomorphic to $S$.
Let $\cH'$ be the union of all these graphs. Then 
$$\lim_{k\to\infty} c^{H_{n_k}}_S= c^{\cH'}_S$$
for any $S$, $|V(S)|\leq K$. Thus the subgraphing $\cH'$ satisfies
the conditions of our lemma. \qed

\vskip 0.2in
\noindent
Now we prove the analogue of Lemma \ref{u3} for graphings.
\begin{lemma} \label{u4}
Let $\hg,\cH,K$ be as above. Then for any $\delta>0$ there exists
$g(\delta)>0$ such that if for a 
hyperfinite graphing $\cG'$, $d_{stat}(\cG',\cG_{\hg})<g(\delta)$ then
$\cG'$ contains a subgraphing $\cH'\subset \cG'$ with components of size at
most $K$ such that
\begin{equation} \label{ronda2}
|c^{\cH'}_S-c^{\cH}_S|<\delta  \,\,\,\mbox{for all}\,\,\, S, |V(S)|\leq K
\end{equation}
\end{lemma}
\proof
Let $d_{stat}(\cG',\cG_{\hg})<f(\delta/2)/2$, where $f$ is the function
in Lemma \ref{u3}. Since $\cG'$ is hyperfinite, it has a subgraphing
$\cK\subset \cG'$ consisting of components of size not greater than some
constant $l>0$. Let us choose a graph sequence $\{Q_n\}^\infty_{n=1}$
such that
\begin{itemize}
\item $\{Q_n\}^\infty_{n=1}$ converges to $\cK$.
\item $d_{stat}(\cK,Q_n)<\frac{f(\delta/2)}{2}$ for all $n\geq 1$.
\end{itemize}
Therefore, $d_{stat}(\cG,Q_n)<f(\delta/2)$ holds for all $n\geq 1$.
Hence by Lemma \ref{u3}, there exist subgraphs $H_n\subset Q_n$ with 
components of size at most $K$ such that
$|c^{H_n}_S-c^{\cH}_S|<\delta/2$ for any $S$, $|V(S)|\leq K$.
By Lemma \ref{becslo}, we have a subgraphing $\cH'\subset \cK$
with components of size at most $K$ satisfying (\ref{ronda2}). \qed

\vskip 0.2in
\noindent
Now we finish the proof of our theorem.
Observe that $\cL(L\cP)\subset [0,1]^{\NN}$ is a compact set.
Call a hyperfinite graphing $\cG$ an $(\e,K)$-graphing if
one can remove an edge set of edge-measure $\e$ to obtain a subgraphing
with components of size at most $K$. By Lemma \ref{u4}, if $\cG\in L\cP$ is
an $(\e,K)$-graphing
then if $d_{stat}(\cG,\cG')$ is small enough then $\cG'$ is 
an $(2\e,K)$-graphing.
So, the theorem follows from compactness. \qed

\vskip 0.2in
\noindent
{\bf Remark:} The reader might ask, whether if $\cP$ is a hyperfinite
family of $(\e,K)$-graphs, then what is the best constant in the Uniformicity
Theorem. As a matter of fact, any constant $\e'>\e$ is good. Indeed,
if $\{Q_n\}^\infty_{n=1}\subset \cP$ is a convergent sequence of 
$(\e,K)$-graphs, then according to the construction in Proposition \ref{masik}
there exists an $(\e',K)$-good limit graphing. So, $\e'$ is a good constant
for the Uniformicity Theorem by Theorem \ref{harmadiktetel}. 
\subsection{The proof of the Equipartition Theorem}
By the Uniformicity Theorem, all elements of $L\cP$ are $(\e,K)$-graphings
for some $K>0$.
Suppose that the theorem does not hold for some $\delta>0$.
Then we have a sequence of graphs $\{G_n\}^\infty_{n=1}, \{H_n\}^\infty_{n=1}
\subset \Gd$ such that $\lim_{n\to\infty} d_{stat}(G_n,H_n)=0$, without
having pairs $\{G_n',H_n'\}^\infty_{n=1}$ satisfying the requirement of
the theorem. Let us pick a convergent graph sequence
$\hg=\{G_{n_k}\}^\infty_{k=1}$. Then $\{H_{n_k}\}^\infty_{k=1}$ tends to
$\cG_{\hg}$ as well. Let $\cH\subset \cG_{\hg}$ be a subgraphing
with components of size at most $K$.
By the Transfer Theorem, we have subgraphs 
$\{G'_{n_k}\subset G_{n_k}\}^\infty_{k=1}$, $\{H'_{n_k}\subset H_{n_k}\}^\infty_{k=1}$
converging to $\cH$. By Lemma \ref{u3}, we can suppose that all the components
of $G'_{n_k}$ and $H'_{n_k}$ have size at most $K$.
Then for large enough $k$,
$$|E(G_{n_k})\backslash E(G'_{n_k})|\leq 2\e |E(G_{n_k})|
\,\,\,\mbox{and} \,\,\, |E(H_{n_k})\backslash E(H'_{n_k})|
\leq 2\e |E(H_{n_k})|$$
Also,
$$\sum_S|c_S^{G_{n_k}'}- c^{\cH}_S|<\frac{\delta}{2}
\,\,\,\mbox{and} \,\,\, 
\sum_S|c_S^{H_{n_k}'}- c^{\cH}_S|<\frac{\delta}{2}$$
leading to a contradiction. \qed
\subsection{The proof of Theorem \ref{newman}}
Let $\e>0$, $\kappa>0$ be constants such that
$(2\e d +\kappa d)<\delta$.
Suppose that $d_{stat}(G,H)<f(\kappa)$, where
$f$ is the function in the Equipartition Theorem. So, we have subgraphs
$G'\subset G$, $H'\subset H$ such that
\begin{itemize}
\item $\sum_S |c^{G'}_S-c^{H'}_S|<\kappa\,.$
\item $|E(G)\backslash E(G')|<2\e |E(G)|\leq \e d n$
\item $|E(H)\backslash E(H')|<2\e |E(H)|\leq \e d n$
\end{itemize}
Then if $c^{G'}_S\leq c^{H'}_S$, we define $\rho: C^{G'}_S\to C^{H'}_S$ to
be a component preserving injective map.
On the other hand, if $c^{G'}_S\geq c^{H'}_S$, then let $D^{G'}_S\subset C^{G'}_S$
be a union of some components such that
$|C^{H'}_S|=|D^{G'}_S|$ and define $\rho:D^{G'}_S\to C^{H'}_S$ to be a
component preserving bijection.
Finally, extend $\rho$ to $V(G)$ arbitrarily.
Observe that
$$|\rho^{-1}(E(H))\triangle E(G)|\leq (2\e d +\kappa d) n\quad \qed $$

\section{Local-global convergence}
The notion of local-global convergence was introduced by Hatami, Lov\'asz
and Szegedy \cite{HLS} (and independently by Bollob\'as and Riordan \cite{BR}
under the name of convergence in the partition metric).

First, let us recall the definition. For $k\geq 2$, let $U^{r,k}_d$
be the finite set of rooted $r$-balls $H$ with
vertex labelings $c:V(H)\to \{1,2,\dots,k\}= [k]$.
Let $G\in \Gd$ be a finite graph.
One can associate to a labeling $c$ a probability distribution
$P_c$ on $\urkd$, where $P_c(\gamma)=p(G,c,\gamma)$, and 
$p(G,c,\gamma)$ is the probability that the $r$-neighborhood of a
random vertex of $G$ is labeled-isomorphic to $\gamma$.
Set 
$$C_k(G):=\cup_{c:V(H)\to [k]}\subset [0,1]^{\urkd}\,.$$
The $k$-th partition pseudodistance of $G$ and $H$ is 
$d_k(G,H):=d_{haus}(C_k(G),C_k(H))$, where $d_{haus}$ is the Hausdorff-distance.
The local-global pseudodistance of $G$ and $H$ is given by
$d_{LG}(G,H)=\sum^\infty_{k=1}\frac{1}{2^{k}} d_k(G,H)\,.$
We can extend the local-global pseudodistance to graphings, as well.
Let $\cG(X,\mu)$ be a graphing of vertex degree bound $d$
and $c:X\to [k]$ be a Borel function.
Then $P_c(\gamma)=\mu(T(\cG,c,\gamma))$, where $(T(\cG,c,\gamma))$
is the set of vertices in $X$ with $r$-neighborhood isomorphic to $\gamma$
(under the labeling induced by $c$).
Let $C_k(\cG)$ be the closure of the set $\cup_c P_c\subset [0,1]^{\urkd}$
and the local-global pseudodistance can be defined as in the case
of finite graphs.
A graph sequence $\{G_n\}^\infty_{n=1}$ converges locally-globally to 
a graphing $\cG$ if for any $k\geq 1$,
$\{C_k(G_n)\}^\infty_{n=1}$ converges to $C_k(\cG)$ in the Hausdorff distance.
Although in general, local-global convergence is much stronger than
the Benjamini-Schramm convergence, for hyperfinite sequences the two
notions coincide (see also Theorem 9.5 \cite{HLS}).
\begin{theorem}\label{localglobal}
If $\{G_n\}^\infty_{n=1}$ is a hyperfinite graph sequence converging
to $\cG$ then it converges to $\cG$ locally-globally.
\end{theorem}
\proof
The following lemma is straightforward and left for the reader.
It states that a small perturbation of a graph is close to the
original graph in the local-global distance.
\begin{lemma} \label{L1}
For any $\e>0$, there exists $\delta>0$ such that
\begin{itemize}
\item If $G\in \Gd$, $H\subset G$, $V(H)=V(G)$ and
$\frac{|E(G\backslash H)|}{|V(G)|}<\delta$ then
$d_{LG}(G,H)<\e$. 
\item If $\cG(X,\mu)$ is a graphing, $\cH\subset \cG$, $V(\cH)\subset V(\cG)$
and $\mu_E(\cG\backslash \cH)<\delta$ then
$d_{LG}(\cG,\cH)<\e$.
\end{itemize}
\end{lemma}
\begin{lemma} \label{L2}
Let $\cH(X,\mu)$ be a graphing with components of size at most $K$.
Let $\{H_n\}^\infty_{n=1}\subset \Gd$ be graphs with
components of size at most $K$ converging to $\cH$. Then
$\lim_{n\to\infty} d_{LG}(H_n,\cH)=0\,.$
\end{lemma}
\proof
By Lemma \ref{u2}, $lim_{n\to\infty} c^{H_n}_S=c^{\cH}_S$ for any
$S$, $|V(S)|\leq K$.
Let $M(K,S,k)$ be the set of all non-isomorphic $k$-labelings of $S$.
A Borel map $c:X\to [k]$ determines
a probability distribution on $M(K,S,k)$, where
$P_c(\beta)=\mu_{\cH}(T(\beta))\,.$ Clearly, for any $\e>0$
there exists $\delta >0$  if
$|c^{H_n}_S-c^{\cH}_S|<\delta$ then we can partition
$C^{H_n}_S$ into parts $C^{H_n}_S=\cup_{\beta\in M(K,S,k)} L(\beta)$
such that $|\frac{L(\beta)}{V(H_n)}-P_c(\beta)|<\e$.
Conversely, for any $\e>0$ there exists $\delta>0$ such
that if $|c_S^{H_n}-c_S^{\cH}|<\delta$ and $c^{H_n}_S=\sum_{\beta\in M(K,S,k)}
l(\beta)$, $l(\beta)\geq 0$, then one can divide $C^{\cH}_S$ into
Borel parts $C^{\cH}_S=\cup_{\beta\in M(K,S,k)} R(\beta)$ in such a way
that each component of $R(\beta)$ is a component of $C^{\cH}_S$
and $|\mu(R(\beta))-l(\beta)|<\e$. Hence the lemma follows, since
the $\e$-ball around $C_k(\cH)$ contains $C_k(H_n)$ and vice versa, the
$\e$-ball around $C_k(H_n)$ contains $C_k(\cH)$. \qed

\vskip 0.2in
\noindent
Now we finish the proof of our theorem. By Lemma \ref{L1} and Lemma \ref{u3},
we have $\cH\subset\cG$ and $\{H_n\subset G_n\}^\infty_{n=1}$
such that
$$d_{LG}(\cG,\cH)<\frac{\e}{3},\,\, d_{LG}(G_n,H_n)<\frac{\e}{3}\,.$$
Hence if $n$ is large, then
$d_{LG}(\cG,G_n)<\e$. \qed

\section{Strong equivalence}
\subsection{The proof of Theorem \ref{harmadiktetel}}
First we define a new pseudo-distance for graphings.
Let $\cG(X,\mu)$, $\cH(Y,\nu)$ be graphings of vertex degree bound $d$.
Then let $d_{strong}(\cG,\cH)$ be the infimum of $\e'$s such that there
exists a measure preserving bijection $\rho:X\to Y$ with
$$\mu_E(\rho^{-1}(E(\cH))\triangle E(\cG))\leq \e\,.$$
So, $\cG$ and $\cH$ is strongly equivalent if $d_{strong}(G,H)=0$.
\begin{lemma} \label{kicsik}
Let $\cH_1(X,\mu)$, $\cH_2(Y,\nu)$ be graphings of degree bound $d$ with
components of size at most $K$. Suppose that
$$\sum_{S, |V(S)|\leq K} |c_S^{\cH_1}-c_S^{\cH_2}|<\kappa\,.$$
Then $d_{strong}(\cH_1, \cH_2)< d\kappa\,.$
\end{lemma}
\proof
Let $S$ be a connected graph of size at most $K$. If
$c^{\cH_1}_S\leq c^{\cH_2}_S$, then define $\rho:C^{\cH_1}_S\to C^{\cH_2}_S$ to be
an injective map preserving the components such that
$\nu(\rho(A))=\mu(A)$ if $A\subset C^{\cH_1}_S$ is a measurable set.
On the other hand, if $c^{\cH_1}_S > c^{\cH_2}_S$, then  let
$D^{\cH_1}_S$ be a Borel set of $X$ such that
\begin{itemize}
\item The components of $D^{\cH_1}_S$ are components of $C^{\cH_1}_S$.
\item $\mu(D^{\cH_1}_S)=c^{\cH_2}_S$.
\end{itemize}
Define $\rho:D_S^{\cH_1}\to C_S^{\cH_2}$ to be a measure-preserving bijection (that
also preserves the components). Then, extend $\rho$ to
a measure-preserving bijection arbitrarily onto the whole space $X$.
Then $\mu_E(\rho^{-1}(E(\cH_2))\triangle E(\cH_1))\leq d\kappa\,.$ \qed

\vskip 0.2in
Now let us finish the proof of Theorem \ref{harmadiktetel}.
Let $Z\subset \cG$ be a set of edges of edge-measure less than $\e/4$, such
that the components of $\cK=\cG\backslash Z$ are of size at most $K$.
Then by the definition of $d_{stat}$ respectively by Lemma \ref{u4}, 
there exists 
some $\delta>0$ such that
\begin{itemize}
\item $d_{stat}(\cG,\cH)<\delta$, then $|\mu_E(\cG)-\nu_E(\cH)|<\e/4\,.$
\item $d_{stat}(\cG,\cH)<\delta$, then $\cH$ contains a subgraphing $\cK'$
such that $\sum_{S, |V(S)|\leq K} |c^{\cK}_S-c^{\cK'}_S|<\frac{\e}{4d}$.
\end{itemize}
By the previous lemma,
\begin{equation}\label{randa3}
d_{strong}(\cK,\cK')<\e/4
\end{equation}
By (\ref{randa3}), $|\mu_E(\cK)-\nu_E(\cK')|<\epsilon/4$.
Thus
$$d_{strong}(\cG,\cH)\leq d_{strong}(\cK,\cK')+\mu_E(\cG\backslash\cK)+
\mu_E(\cH\backslash\cK')\leq \frac{\e}{4} + \frac{\e}{4} + 
\frac{\e}{2}\leq\e\quad\qed$$
\subsection{Rokhlin Lemma for non-free actions}
Let $\Gamma$ be a finitely generated amenable group with
a symmetric generating system.
Ornstein and Weiss \cite{OW} proved the
following version of the classical Rokhlin Lemma.
If $\Gamma\curvearrowright (X,\mu)$, $\Gamma\curvearrowright(Y,\nu)$
are two probability measure preserving essentially free actions, then
they are strongly equivalent. That is for any $\e>0$  there exists
a measure preserving bijection
$\rho_\e:X\to Y$ such that
$$\mu(\{x\in X,\,\mid\, \rho_\e(sx)=s\rho_\e(x) \,\,\mbox{for any $s\in S$}
\,\,\})>1-\e\,.$$
The goal of this subsection is to show how one can deduce
the general (non-free) version of the statement above using
Theorem \ref{harmadiktetel}.
First, let us recall the notion of the {\it type} of an action
(\cite{AE}, \cite{TD}).
Let $\FF_n$ be the free group on $n$-generators
$\{s_1, s_2,\dots, s_n\}$. Let $\alpha=\FF_n\curvearrowright(X,\mu)$
be a not necessarily free action of $\FF_n$.
Note that any free action of an $n$-element generated group $\Gamma$
can be viewed as a non-free action of $\FF_n$.
Let $\Sigma_n$ be the space of all rooted Schreier graphs of
transitive actions of $\FF_n$ on countable sets. Note that
the elements of $\Sigma_n$ are connected rooted graphs
with edge labels from $\{s_1,s_2,\dots, s_n,s^{-1}_1, s^{-1}_2,\dots, s^{-1}_n\}$
where the edge $(x,s_ix)$ is labeled by $s_i$.
The space $\Sigma_n$ is compact and $\FF_n$ acts on $\Sigma_n$ 
continuously by changing the roots. Following \cite{AGV},
we call the $\FF_n$-invariant measures on $\Sigma_n$
{\it invariant random subgroups} (IRS). Let
$\alpha:\FF_n\curvearrowright(X,\mu)$ be a p.m.p. Borel action.
The type of $\alpha$ is an IRS defined the following way.
Let $\pi_\alpha:X\to\Sigma_n$ be the
map that maps $x\in X$ to the Schreier graph of its orbit (with root $x$).
The type of $\alpha$, $type(\alpha)$ is the invariant
measure $(\pi_\alpha)_\star(\mu)$.
Now we state the non-free version of the amenable Rokhlin Lemma.
Note that a version (stably weak equivalence of the actions)
of the result is proved in \cite[Theorem 1.8]{TD}.
\begin{theorem}
\label{rohlin}
If $\alpha,\beta: \FF_n\curvearrowright (X,\mu)$ are hyperfinite
actions (the underlying graphings are hyperfinite) and
$type(\alpha)=type(\beta)$ then $\alpha$ and $\beta$ are strongly
equivalent.
\end{theorem}
\proof The idea of the proof is that for each action $\alpha$
we construct an (unlabeled) graphing $\cG_\alpha$ such that
$type(\alpha)=type(\beta)$ if and only if $d_{stat}(\cG_\alpha,\cH_\alpha)$.
One should note that if the orbits have no rooted automorphisms, then
the graphing of $\alpha$ would fit for this purpose. Again, we only need to
handle the symmetries. 
First, let $\cG^{\alpha}(X,\mu)$ be the graphing of our action.
We will ``add'' marker graphs to $\cG^{\alpha}$ in order to encode
the action.
The marker graph for $s_i$ is a path $P_i$ of path-length $i$ (that is
of $i+1$-vertices). The additional marker graph for a vertex in $X$ is
the path $P_{n+1}$. The construction of $\cG_\alpha$ goes as follows.

\vskip 0.2in
\noindent
{\bf Step 1.} Stick a graph $P_{n+1}$ to each vertex of $x\in X$ (the vertices
of $X$ will be called ``original'' vertices).
This means that we identify an endpoint of $P_{n+1}$ with $x$.
In this way, we obtain a new graphing $\cG^\alpha_1(X_1,\mu_1)$.
Here $X_1$ is the union of $n+2$-copies of $X$. We normalize $\mu_1$ in order
to get a probability measure.

\vskip 0.2in
\noindent
{\bf Step 2.} Now we divide each edge $(x,s_i x)$ of the original graphing
$\cG^\alpha$ into three parts by adding two new vertices. In this way,
we obtain the graphing $\cG_2$ from $\cG_1$.
Note that if $x=s_ix$ we do not make any subdivison
(we do not consider loops). Also, if $s_ix=s_jx$ then the
edges $(x,s_ix)$ and $(x,s_jx)$ coincide.

\vskip 0.2in
\noindent
{\bf Step 3.} In the final step we encode the action.
For each $1\leq i \leq n$ we stick a marker graph $P_i$
to the vertex next to $x$ on the path $x,s_ix$, where $x$ is
an original vertex.
The resulting graphing is $\cG_\alpha(X_\alpha,\mu_\alpha)$ 
(the fact that it is measure-preserving
Borel graph follows immediately from the invariance of the action $\alpha$).
By looking at the $3n$-ball around a vertex of $X_\alpha$ we can see whether
it is an original vertex or not. In fact by looking at the
$3nr$-ball around such a vertex we can reconstruct the labeled $r$-ball
of the original labeled graphing $\cG^\alpha$.
It is not hard to see that $type(\alpha)=type(\beta)$ if and only if
$d_{stat}(\cG_\alpha,\cG_\beta)=0$.
Hence if $type(\alpha)=type(\beta)$, by Theorem \ref{harmadiktetel},
$\cG_\alpha$ is strongly equivalent to $\cG_\beta$. This implies
the strong equivalence of the actions $\alpha$ and $\beta$. \qed.

\vskip 1in
\noindent
Alfred Renyi Institute of the Hungarian Academy of Sciences

\noindent
and

\noindent
\'Ecole Polytechnique F\'ed\'erale de Lausanne, EPFL

\vskip 0.3in \noindent
elek@renyi.hu

\end{document}